\documentclass{article}%
\usepackage{amsfonts}
\usepackage{amsmath}
\usepackage{amssymb}
\usepackage{graphicx}%
\newtheorem{theorem}{Theorem}

\newtheorem{condition}[theorem]{Condition}

\newtheorem{corollary}[theorem]{Corollary}

\newtheorem{definition}[theorem]{Definition}
\newtheorem{example}[theorem]{Example}

\newtheorem{proposition}[theorem]{Proposition}
\newtheorem{remark}[theorem]{Remark}

\newenvironment{proof}[1][Proof]{\noindent\textbf{#1.} }{\ \rule{0.5em}{0.5em}}
\begin{document}

\title{On a factorization of second order elliptic operators and applications}
\author{Vladislav V. Kravchenko\\Secci\'{o}n de Posgrado e Investigaci\'{o}n\\Escuela Superior de Ingenier\'{\i}a Mec\'{a}nica y El\'{e}ctrica\\Instituto Polit\'{e}cnico Nacional\\C.P.07738 M\'{e}xico D.F., \\MEXICO\\e-mail: vkravchenko@ipn.mx}
\maketitle

\begin{abstract}
We show that given a nonvanishing particular solution of the equation
\begin{equation}
(\operatorname{div}p\operatorname{grad}+q)u=0, \label{1}%
\end{equation}
the corresponding differential operator can be factorized into a product of
two first order operators. The factorization allows us to reduce the equation
(1) to a first order equation which in a two-dimensional case is the Vekua
equation of a special form. Under quite general conditions on the coefficients
$p$ and $q$ we obtain an algorithm which allows us to construct in explicit
form the positive formal powers (solutions of the Vekua equation generalizing
the usual powers $(z-z_{0})^{n}$, $n=0,1,\ldots$). This result means that
under quite general conditions one can construct an infinite system of exact
solutions of (1) explicitly, and moreover, at least when $p$ and $q$ are real
valued this system will be complete in $\ker(\operatorname{div}%
p\operatorname{grad}+q)$ in the sense that any solution of (1) in a simply
connected domain $\Omega$ can be represented as an infinite series of obtained
exact solutions which converges uniformly on any compact subset of $\Omega$.

Finally we give a similar factorization of the operator $(\operatorname{div}%
p\operatorname{grad}+q)$ in a multidimensional case and obtain a natural
generalization of the Vekua equation which is related to second order
operators in a similar way as its two-dimensional prototype does.

\end{abstract}

\section{Introduction}

Consider the one-dimensional stationary Schr\"{o}dinger equation%
\begin{equation}
y^{\prime\prime}+v(x)y=0. \label{Schrod1}%
\end{equation}
It is well known that given a nonvanishing particular solution $y_{0}$ of
(\ref{Schrod1}), the Schr\"{o}dinger operator can be factorized%
\begin{equation}
\partial^{2}+v(x)=\left(  \partial+\frac{y_{0}^{\prime}(x)}{y_{0}(x)}\right)
\left(  \partial-\frac{y_{0}^{\prime}(x)}{y_{0}(x)}\right)
\label{FactSchrod1}%
\end{equation}
and as a consequence the general solution of (\ref{Schrod1}) can be obtained.

In \cite{Krpseudoan} it was shown that a similar factorization of a stationary
Schr\"{o}dinger operator is available in two dimensions. One of the main
results of the present work is the factorization of a more general
two-dimensional elliptic operator $(\operatorname{div}p\operatorname{grad}+q)$
in a form similar to (\ref{FactSchrod1}). The factorization allows us to
reduce the equation
\begin{equation}
(\operatorname{div}p\operatorname{grad}+q)u=0 \label{maineqintro}%
\end{equation}
to a Vekua equation of a special form. This Vekua equation becomes bicomplex
if $p$ or $q$ are complex functions. One complex component of a solution of
the Vekua equation (its real part when $p$ and $q$ are real valued)
necessarily satisfies equation (\ref{maineqintro}) and the other satisfies an
associated equation having the form of (\ref{maineqintro}) but with different
coefficients $p$ and $q$. This situation is a generalization of the fact that
real and imaginary parts of an analytic function are harmonic, and likewise
for any harmonic function in a simply connected domain its harmonic conjugate
can be constructed, we obtain explicit formulas for constructing
\textquotedblleft conjugate\textquotedblright\ solutions of the associated
equations of the form (\ref{maineqintro}). In the case $p\equiv1$ and
$q\equiv0$ these formulas turn into the well known from complex analysis
formulas for construction of conjugate harmonic functions.

In \cite{Krpseudoan} we established that under quite general conditions the
positive formal powers corresponding to the Vekua equation (the pseudoanalytic
functions generalizing the usual powers $(z-z_{0})^{n}$, $n=0,1,\ldots$) can
be constructed explicitly. Here we develop this result and obtain that by a
given nonvanishing particular solution of equation (\ref{maineqintro}) under
quite general conditions one can construct an infinite system of exact
solutions of (\ref{maineqintro}) explicitly, and moreover, at least when $p$
and $q$ are real valued this system will be complete in $\ker
(\operatorname{div}p\operatorname{grad}+q)$ in the sense that any solution of
(\ref{maineqintro}) in a simply connected domain $\Omega$ can be represented
as an infinite series of obtained exact solutions which converges uniformly on
any compact subset of $\Omega$.

In the final part of the present work we obtain a factorization of the
operator $(\operatorname{div}p\operatorname{grad}+q)$ in a multidimensional
situation and reduce (\ref{maineqintro}) to a first order equation which
generalizes the Vekua equation.

\section{Preliminaries}

A bicomplex number has the form
\[
q=Q_{1}+Q_{2}k
\]
where $Q_{1}=q_{0}+iq_{1}$, $Q_{2}=q_{2}+iq_{3}$, $q_{j}\in\mathbb{R}$,
$j=\overline{0,3}$; $i^{2}=k^{2}=-1$, $ik=ki$.

We will say that $Q_{1}$ and $Q_{2}$ are complex components of the bicomplex
number $q$. Denote $\overline{q}=Q_{1}-Q_{2}k$. The corresponding conjugation
operator we denote by $C$: $Cq=\overline{q}$. We will say that $q$ is scalar
if $Q_{2}=0$. That is in this case $q$ is a usual complex number. $Q_{2}$ is
called the vector part of $Q$. Sometimes we use the notation: $Q_{1}%
=\operatorname*{Sc}(Q)$, $Q_{2}=\operatorname*{Vec}(Q)$.

The set of bicomplex zero divisors, that is of nonzero elements $q=Q_{1}%
+Q_{2}k$, $\left\{  Q_{1},Q_{2}\right\}  \subset\mathbb{C}$ such that
\begin{equation}
q\overline{q}=\left(  Q_{1}+Q_{2}k\right)  \left(  Q_{1}-Q_{2}k\right)  =0
\label{Defzerodiv}%
\end{equation}
we denote by $\mathfrak{S}$.

We will consider the variable $z=x+ky$, where $x$ and $y$ are real variables
and the corresponding differential operators%
\[
\partial_{\overline{z}}=\frac{1}{2}(\partial_{x}+k\partial_{y})\text{\quad
and\quad}\partial_{z}=\frac{1}{2}(\partial_{x}-k\partial_{y}).
\]
Notation $W_{\overline{z}}$ or $W_{z}$ means the application of $\partial
_{\overline{z}}$ or $\partial_{z}$ respectively to a bicomplex function $W(z)$.

Note that if we consider the equation
\begin{equation}
\partial_{z}\varphi=\Phi\label{gradient}%
\end{equation}
in a whole complex plane or in a convex domain, where $\Phi=\Phi_{1}+k\Phi
_{2}$ is a given bicomplex valued function such that its scalar part $\Phi
_{1}$ and vector part $\Phi_{2}$ satisfy the equation
\begin{equation}
\partial_{y}\Phi_{1}+\partial_{x}\Phi_{2}=0, \label{casirot}%
\end{equation}
then there exists a scalar solution of (\ref{gradient}) which can be
reconstructed up to an arbitrary scalar constant $c$ in the following way%
\begin{equation}
\varphi(x,y)=2\left(  \int_{x_{0}}^{x}\Phi_{1}(\eta,y)d\eta-\int_{y_{0}}%
^{y}\Phi_{2}(x_{0},\xi)d\xi\right)  +c \label{Antigr}%
\end{equation}
where $(x_{0},y_{0})$ is an arbitrary fixed point in the domain of interest.

By $A$ we denote the integral operator in (\ref{Antigr}):%
\[
A[\Phi](x,y)=2\left(  \int_{x_{0}}^{x}\Phi_{1}(\eta,y)d\eta-\int_{y_{0}}%
^{y}\Phi_{2}(x_{0},\xi)d\xi\right)  +c.
\]
Note that formula (\ref{Antigr}) can be easily extended to any simply
connected domain by considering the integral along an arbitrary rectifiable
curve $\Gamma$ leading from $(x_{0},y_{0})$ to $(x,y)$%
\[
\varphi(x,y)=2\left(  \int_{\Gamma}\Phi_{1}dx-\Phi_{2}dy\right)  +c.
\]
Thus if $\Phi$ satisfies (\ref{casirot}), there exists a family of scalar
functions $\varphi$ such that $\partial_{z}\varphi=\Phi$, given by the formula
$\varphi=A[\Phi]$.

In a similar way we define the operator $\overline{A}$ corresponding to
$\partial_{\overline{z}}$:%
\[
\overline{A}[\Phi](x,y)=2\left(  \int_{x_{0}}^{x}\Phi_{1}(\eta,y)d\eta
+\int_{y_{0}}^{y}\Phi_{2}(x_{0},\xi)d\xi\right)  +c.
\]

\section{Solutions of second order elliptic equations as scalar parts of
bicomplex pseudoanalytic functions\label{SectSchr}}

Consider the equation%
\begin{equation}
\left(  -\Delta+\nu\right)  f=0 \label{Schrod}%
\end{equation}
in some domain $\Omega\subset\mathbf{R}^{2}$, where $\Delta=\frac{\partial
^{2}}{\partial x^{2}}+\frac{\partial^{2}}{\partial y^{2}}$, $\nu$ and $f$ are
complex valued (in our terms scalar) functions. We assume that $f$ is a twice
continuously differentiable function.

\begin{theorem}
\label{Thfact}Let $f$ be a nonvanishing in $\Omega$ particular solution of
(\ref{Schrod}). Then for any complex valued (scalar) function $\varphi\in
C^{2}(\Omega)$ the following equalities hold%
\begin{equation}
\frac{1}{4}\left(  \Delta-\nu\right)  \varphi=\left(  \partial_{\overline{z}%
}+\frac{f_{z}}{f}C\right)  \left(  \partial_{z}-\frac{f_{z}}{f}C\right)
\varphi=\left(  \partial_{z}+\frac{f_{\overline{z}}}{f}C\right)  \left(
\partial_{\overline{z}}-\frac{f_{\overline{z}}}{f}C\right)  \varphi.
\label{fact}%
\end{equation}

\end{theorem}

\begin{proof}
Consider%
\begin{align*}
\left(  \partial_{\overline{z}}+\frac{f_{z}}{f}C\right)  \left(  \partial
_{z}-\frac{f_{z}}{f}C\right)  \varphi &  =\frac{1}{4}\Delta\varphi
-\frac{\left\vert \partial_{z}f\right\vert ^{2}}{f^{2}}\varphi-\partial
_{\overline{z}}\left(  \frac{\partial_{z}f}{f}\right)  \varphi\\
&  =\frac{1}{4}(\Delta\varphi-\frac{\Delta f}{f}\varphi)=\frac{1}{4}\left(
\Delta-\nu\right)  \varphi.
\end{align*}

Application of $C$ gives us the last part of (\ref{fact}).
\end{proof}

In the case of a real valued potential $\nu$ this theorem was proved in
\cite{Krpseudoan}.

The following statement is known in a form of a substitution (see, e.g.,
\cite{Nachman}). Here we formulate it as an operator relation.

\begin{proposition}
\label{PropSchr}Let $p$ and $q$ be complex valued functions, $p\in
C^{2}(\Omega)$ and $p\neq0$ in $\Omega$. Then
\begin{equation}
\operatorname{div}p\operatorname{grad}+q=p^{1/2}(\Delta-r)p^{1/2}\text{\qquad
in }\Omega, \label{firstfact}%
\end{equation}
where
\[
r=\frac{\Delta p^{1/2}}{p^{1/2}}-\frac{q}{p}.
\]

\end{proposition}

\begin{proof}
The easily verified relation
\begin{equation}
\operatorname{div}p\operatorname{grad}=p^{1/2}(\Delta-\frac{\Delta p^{1/2}%
}{p^{1/2}})p^{1/2} \label{factUhlm}%
\end{equation}
is well known (see, e.g., \cite{UhlmannDevelopments}). Adding to both sides of
(\ref{factUhlm}) the term $q$ (and representing it on the right-hand side as
$p^{1/2}\left(  q/p\right)  p^{1/2}$) gives us (\ref{firstfact}).
\end{proof}

\begin{theorem}
\label{ThFactGenSchr}Let $u_{0}$ be a nonvanishing in $\Omega$ particular
solution of the equation
\begin{equation}
(\operatorname{div}p\operatorname{grad}+q)u=0\text{\qquad in }\Omega\text{.}
\label{maineq}%
\end{equation}
Then under the conditions of proposition \ref{PropSchr} for any complex valued
(scalar) continuously twice differentiable function $\varphi$ the following
equality holds%
\begin{equation}
\frac{1}{4}(\operatorname{div}p\operatorname{grad}+q)\varphi=p^{1/2}\left(
\partial_{z}+\frac{f_{\overline{z}}}{f}C\right)  \left(  \partial
_{\overline{z}}-\frac{f_{\overline{z}}}{f}C\right)  p^{1/2}\varphi,
\label{mainfact}%
\end{equation}
where
\begin{equation}
f=p^{1/2}u_{0}. \label{fandu}%
\end{equation}

\end{theorem}

\begin{proof}
This is based on (\ref{fact}). From (\ref{firstfact}) we have that if $u_{0}$
is a solution of (\ref{maineq}) then the function (\ref{fandu}) is a solution
of the equation%
\begin{equation}
(\Delta-r)f=0. \label{Schrodr}%
\end{equation}
Then combining (\ref{firstfact}) and (\ref{fact}) we obtain (\ref{mainfact}).
\end{proof}

\begin{remark}
According to (\ref{factUhlm}), $\Delta-r=f^{-1}\operatorname{div}%
f^{2}\operatorname{grad}f^{-1}$ where $f$ is a solution of (\ref{Schrodr}).
Then from (\ref{firstfact}) we have%
\begin{equation}
\operatorname{div}p\operatorname{grad}+q=p^{1/2}f^{-1}\operatorname{div}%
f^{2}\operatorname{grad}f^{-1}p^{1/2}. \label{reldivgrad}%
\end{equation}
Taking into account (\ref{fandu}) we obtain
\[
\operatorname{div}p\operatorname{grad}+q=u_{0}^{-1}\operatorname{div}%
pu_{0}^{2}\operatorname{grad}u_{0}^{-1}\text{\qquad in }\Omega\text{.}%
\]

\end{remark}

\begin{remark}
Let $q\equiv0$. Then $u_{0}$ can be chosen as $u_{0}\equiv1$. Hence
(\ref{mainfact}) gives us the equality
\[
\frac{1}{4}\operatorname{div}(p\operatorname{grad}\varphi)=p^{1/2}\left(
\partial_{z}+\frac{\partial_{\overline{z}}p^{1/2}}{p^{1/2}}C\right)  \left(
\partial_{\overline{z}}-\frac{\partial_{\overline{z}}p^{1/2}}{p^{1/2}%
}C\right)  (p^{1/2}\varphi).
\]

\end{remark}

In what follows we suppose that in $\Omega$ there exists a nonvanishing
particular solution of (\ref{maineq}) which we denote by $u_{0}$.

Let $f$ be a scalar function of $x$ and $y$. Consider the bicomplex Vekua
equation
\begin{equation}
W_{\overline{z}}=\frac{f_{\overline{z}}}{f}\overline{W}\text{\qquad in }%
\Omega\text{.} \label{Vekuamain}%
\end{equation}

Denote $W_{1}=\operatorname{Sc}W$ and $W_{2}=\operatorname{Vec}W.$

\begin{remark}
\cite{Krpseudoan} \label{RemAnotherForm} Equation (\ref{Vekuamain}) can be
written as follows%
\begin{equation}
f\partial_{\overline{z}}(f^{-1}W_{1})+kf^{-1}\partial_{\overline{z}}%
(fW_{2})=0. \label{anotherformVekua}%
\end{equation}

\end{remark}

\begin{theorem}
\label{ThConjugate} Let $W=W_{1}+W_{2}k$ be a solution of (\ref{Vekuamain}).
Then $U=f^{-1}W_{1}$ is a solution of the equation
\begin{equation}
\operatorname*{div}(f^{2}\nabla U)=0\qquad\text{in }\Omega, \label{divf2}%
\end{equation}
and $V=fW_{2}$ is a solution of the equation
\begin{equation}
\operatorname*{div}(f^{-2}\nabla V)=0\qquad\text{in }\Omega, \label{divf-2}%
\end{equation}
the function $W_{1}$ is a solution of the stationary Schr\"{o}dinger equation
\begin{equation}
-\Delta W_{1}+r_{1}W_{1}=0\qquad\text{in }\Omega\label{Schr1}%
\end{equation}
with $r_{1}=\Delta f/f,$ and $W_{2}$ is a solution of the associated
Schr\"{o}dinger equation
\begin{equation}
-\Delta W_{2}+r_{2}W_{2}=0\qquad\text{in }\Omega\label{Schr2}%
\end{equation}
where $r_{2}=2(\nabla f)^{2}/f^{2}-r_{1}$ and $(\nabla f)^{2}=f_{x}^{2}%
+f_{y}^{2}$.
\end{theorem}

\begin{proof}
To prove the first part of the theorem we use the form of equation
(\ref{Vekuamain}) given in Remark \ref{RemAnotherForm}. Multiplying
(\ref{anotherformVekua}) by $f$ and applying $\partial_{z}$ gives%
\[
\partial_{z}f^{2}\partial_{\overline{z}}\left(  f^{-1}W_{1}\right)  +\frac
{k}{4}\Delta\left(  fW_{2}\right)  =0
\]
from where we have that $\operatorname{Sc}\left(  \partial_{z}f^{2}%
\partial_{\overline{z}}\left(  f^{-1}W_{1}\right)  \right)  =0$ which is
equivalent to (\ref{divf2}) where $U=f^{-1}W_{1}$.

Multiplying (\ref{anotherformVekua}) by $f^{-1}$ and applying $\partial_{z}$
gives%
\[
\frac{1}{4}\Delta\left(  f^{-1}W_{1}\right)  +k\partial_{z}f^{-2}%
\partial_{\overline{z}}\left(  fW_{2}\right)  =0
\]
from where we have that $\operatorname{Sc}\left(  \partial_{z}f^{-2}%
\partial_{\overline{z}}\left(  fW_{2}\right)  \right)  =0$ which is equivalent
to (\ref{divf-2}) where $V=fW_{2}$.

From (\ref{factUhlm}) we have
\[
\left(  \Delta-r_{1}\right)  W_{1}=f^{-1}\operatorname*{div}(f^{2}%
\nabla\left(  f^{-1}W_{1}\right)  ).
\]
Hence from the just proven equation (\ref{divf2}) we obtain that $W_{1}$ is a
solution of (\ref{Schr1}).

In order to obtain equation (\ref{Schr2}) for $W_{2}$ it should be noted that
\[
f\operatorname*{div}(f^{-2}\nabla(fW_{2}))=\left(  \Delta-r_{2}\right)
W_{2}.
\]

\end{proof}

In the case of a real valued function $f$ the relation between solutions of
(\ref{Vekuamain}) and equations (\ref{divf2}), (\ref{divf-2}) was observed in
\cite{KrOviedo06}, and between solutions of (\ref{Vekuamain}) and equations
(\ref{Schr1}), (\ref{Schr2}) in \cite{Krpseudoan}.

\begin{remark}
\label{RelAstalaPaivarinta} Observe that the pair of functions
\begin{equation}
F=f\quad\text{and\quad}G=\frac{k}{f} \label{genpair}%
\end{equation}
is a generating pair for (\ref{Vekuamain}). This allows us to rewrite
(\ref{Vekuamain}) in the form of an equation for pseudoanalytic functions of
second kind%
\begin{equation}
\varphi_{\overline{z}}f+\psi_{\overline{z}}\frac{k}{f}=0,
\label{Vekuamain2kind}%
\end{equation}
where $\varphi$ and $\psi$ are scalar functions. If $\varphi$ and $\psi$
satisfy (\ref{Vekuamain2kind}) then $W=\varphi f+\psi\frac{k}{f}$ is a
solution of (\ref{Vekuamain}) and vice versa.

Denote $w=\varphi+\psi k$. Then from (\ref{Vekuamain2kind}) we have%
\[
(w+\overline{w})_{\overline{z}}f+(w-\overline{w})_{\overline{z}}\frac{1}%
{f}=0,
\]
which is equivalent to the equation
\begin{equation}
w_{\overline{z}}=\frac{1-f^{2}}{1+f^{2}}\overline{w}_{\overline{z}}
\label{APform}%
\end{equation}
The relation between (\ref{APform}) and (\ref{divf2}), (\ref{divf-2}) in the
case of a real valued function $f^{2}$ was observed in \cite{AstalaPaivarinta}
and resulted to be essential for solving the Calder\'{o}n problem in the plane.
\end{remark}

\begin{theorem}
Let $W=W_{1}+W_{2}k$ be a solution of (\ref{Vekuamain}). Assume that
$f=p^{1/2}u_{0}$, where $u_{0}$ is a nonvanishing solution of (\ref{maineq})
in $\Omega$. Then $u=p^{-1/2}W_{1}$ is a solution of (\ref{maineq}) in
$\Omega$, and $v=p^{1/2}W_{2}$ is a solution of the equation
\begin{equation}
(\operatorname*{div}\frac{1}{p}\operatorname*{grad}+q_{1})v=0\qquad\text{in
}\Omega, \label{assocmaineq}%
\end{equation}
where
\begin{equation}
q_{1}=-\frac{1}{p}\left(  \frac{q}{p}+2\left\langle \frac{\nabla p}{p}%
,\frac{\nabla u_{0}}{u_{0}}\right\rangle +2\left(  \frac{\nabla u_{0}}{u_{0}%
}\right)  ^{2}\right)  . \label{q1}%
\end{equation}

\end{theorem}

\begin{proof}
According to theorem \ref{ThConjugate}, the function $f^{-1}W_{1}$ is a
solution of (\ref{divf2}). From (\ref{reldivgrad}) we have that
\[
p^{-1/2}\left(  \operatorname{div}p\operatorname{grad}+q\right)
(p^{-1/2}W_{1})=f^{-1}\operatorname*{div}(f^{2}\nabla\left(  f^{-1}%
W_{1}\right)  )
\]
from which we obtain that $u=p^{-1/2}W_{1}$ is a solution of (\ref{maineq}).

In order to obtain the second assertion of the theorem, let us show that%
\[
p^{1/2}(\operatorname*{div}\frac{1}{p}\operatorname*{grad}+q_{1}%
)(p^{1/2}\varphi)=f\operatorname*{div}(f^{-2}\nabla(f\varphi))
\]
for any scalar $\varphi\in C^{2}(\Omega)$. According to (\ref{factUhlm}),%
\[
f\operatorname*{div}(f^{-2}\nabla(f\varphi))=\left(  \Delta-\frac{\Delta
f^{-1}}{f^{-1}}\right)  \varphi=\left(  \Delta-r_{2}\right)  \varphi.
\]
Straightforward calculation gives us the following equality%
\[
\frac{\Delta f^{-1}}{f^{-1}}=\frac{3}{4}\left(  \frac{\nabla p}{p}\right)
^{2}-\frac{1}{2}\frac{\Delta p}{p}+\left\langle \frac{\nabla p}{p}%
,\frac{\nabla u_{0}}{u_{0}}\right\rangle -\frac{\Delta u_{0}}{u_{0}}+2\left(
\frac{\nabla u_{0}}{u_{0}}\right)  ^{2}.
\]
From the condition that $u_{0}$ is a solution of (\ref{maineq}) we obtain the
equality%
\[
-\frac{\Delta u_{0}}{u_{0}}=\frac{q}{p}+\left\langle \frac{\nabla p}{p}%
,\frac{\nabla u_{0}}{u_{0}}\right\rangle .
\]
Thus,
\[
\frac{\Delta f^{-1}}{f^{-1}}=\frac{3}{4}\left(  \frac{\nabla p}{p}\right)
^{2}-\frac{1}{2}\frac{\Delta p}{p}+2\left\langle \frac{\nabla p}{p}%
,\frac{\nabla u_{0}}{u_{0}}\right\rangle +\frac{q}{p}+2\left(  \frac{\nabla
u_{0}}{u_{0}}\right)  ^{2}.
\]
Notice that
\[
\frac{\Delta p^{-1/2}}{p^{-1/2}}=\frac{3}{4}\left(  \frac{\nabla p}{p}\right)
^{2}-\frac{1}{2}\frac{\Delta p}{p}.
\]
Then
\[
\frac{\Delta f^{-1}}{f^{-1}}=\frac{\Delta p^{-1/2}}{p^{-1/2}}+2\left\langle
\frac{\nabla p}{p},\frac{\nabla u_{0}}{u_{0}}\right\rangle +\frac{q}%
{p}+2\left(  \frac{\nabla u_{0}}{u_{0}}\right)  ^{2}.
\]
Now taking $q_{1}$ in the form (\ref{q1}) we obtain the result from
(\ref{firstfact}).
\end{proof}

\begin{theorem}
\label{PrTransform}\cite{Krpseudoan} Let $W_{1}$ be a solution of
(\ref{Schr1}) in a simply connected domain $\Omega$. Then the function
$W_{2},$ solution of (\ref{Schr2}) such that $W=W_{1}+W_{2}k$ is a solution of
(\ref{Vekuamain}), is constructed according to the formula%
\begin{equation}
W_{2}=f^{-1}\overline{A}(kf^{2}\partial_{\overline{z}}(f^{-1}W_{1})).
\label{transfDarboux}%
\end{equation}
It is unique up to an additive term $cf^{-1}$ where $c$ is an arbitrary
complex constant.

Given a solution $W_{2}$ of (\ref{Schr2}), the corresponding solution $W_{1}$
of (\ref{Schr1}) such that $W=W_{1}+W_{2}k$ is a solution of (\ref{Vekuamain}%
), is constructed as follows%
\begin{equation}
W_{1}=-f\overline{A}(kf^{-2}\partial_{\overline{z}}(fW_{2}))
\label{transfDarbouxinv}%
\end{equation}
up to an additive term $cf.$
\end{theorem}

\begin{remark}
When $\nu\equiv0$ and $f_{0}\equiv1$, equalities (\ref{transfDarboux}) and
(\ref{transfDarbouxinv}) turn into the well known in complex analysis formulas
for constructing conjugate harmonic functions.
\end{remark}

\begin{corollary}
\label{ThDarbouxstatic} Let $U$ be a solution of (\ref{divf2}). Then a
solution $V$ of (\ref{divf-2}) such that
\[
W=fU+kf^{-1}V
\]
is a solution of (\ref{Vekuamain}), is constructed according to the formula%
\[
V=\overline{A}(kf^{2}U_{\overline{z}}).
\]
It is unique up to an additive complex constant. Conversely, given a solution
$V$ of (\ref{divf-2}), the corresponding solution $U$ of (\ref{divf2}) can be
constructed as follows:%
\[
U=-\overline{A}(kf^{-2}V_{\overline{z}}).
\]
It is unique up to an additive complex constant.
\end{corollary}

\begin{proof}
Consists in substitution of $W_{1}=fU$ and of $W_{2}=f^{-1}V$ into
(\ref{transfDarboux}) and (\ref{transfDarbouxinv}).\bigskip
\end{proof}

\begin{corollary}
\label{CorConjugate}Let $f=p^{1/2}u_{0}$, where $u_{0}$ is a nonvanishing
solution of (\ref{maineq}) in a simply connected domain $\Omega$ and $u$ be a
solution of (\ref{maineq}). Then a solution $v$ of (\ref{assocmaineq}) with
$q_{1}$ defined by (\ref{q1}) such that $W=p^{1/2}u+kp^{-1/2}v$ is a solution
of (\ref{Vekuamain}), is constructed according to the formula%
\[
v=u_{0}^{-1}\overline{A}(kpu_{0}^{2}\partial_{\overline{z}}(u_{0}^{-1}u)).
\]
It is unique up to an additive term $cu_{0}^{-1}$ where $c$ is an arbitrary
complex constant.

Let $v$ be a solution of (\ref{assocmaineq}), then the corresponding solution
$u$ of (\ref{maineq}) such that $W=p^{1/2}u+kp^{-1/2}v$ is a solution of
(\ref{Vekuamain}), is constructed according to the formula%
\[
u=-u_{0}\overline{A}(kp^{-1}u_{0}^{-2}\partial_{\overline{z}}(u_{0}v)).
\]

\end{corollary}

\begin{proof}
Consists in substitution of $f=p^{1/2}u_{0}$, $W_{1}=p^{1/2}u$ and
$W_{2}=p^{-1/2}v$ into (\ref{transfDarboux}) and (\ref{transfDarbouxinv}%
).\bigskip
\end{proof}

\section{Some definitions and results from pseudoanalytic theory for bicomplex
functions}

\subsection{Generating pair, derivative and antiderivative}

Following \cite{Berskniga} we introduce the notion of a bicomplex generating pair.

\begin{definition}
A pair of bicomplex functions $F=F_{1}+F_{2}k$ and $G=G_{1}+G_{2}k$,
possessing in $\Omega$ partial derivatives with respect to the real variables
$x$ and $y$ is said to be a generating pair if it satisfies the inequality%
\[
\operatorname{Vec}(\overline{F}G)\neq0\qquad\text{in }\Omega.
\]
The following expressions are called characteristic coefficients of the pair
$(F,G)$
\[
a_{(F,G)}=-\frac{\overline{F}G_{\overline{z}}-F_{\overline{z}}\overline{G}%
}{F\overline{G}-\overline{F}G},\qquad b_{(F,G)}=\frac{FG_{\overline{z}%
}-F_{\overline{z}}G}{F\overline{G}-\overline{F}G},
\]

\end{definition}

\[
A_{(F,G)}=-\frac{\overline{F}G_{z}-F_{z}\overline{G}}{F\overline{G}%
-\overline{F}G},\qquad B_{(F,G)}=\frac{FG_{z}-F_{z}G}{F\overline{G}%
-\overline{F}G}.
\]

Every bicomplex function $W$ defined in a subdomain of $\Omega$ admits the
unique representation $W=\phi F+\psi G$ where the functions $\phi$ and $\psi$
are scalar.

The $(F,G)$-derivative $\overset{\cdot}{W}=\frac{d_{(F,G)}W}{dz}$ of a
function $W$ exists and has the form
\begin{equation}
\overset{\cdot}{W}=\phi_{z}F+\psi_{z}G=W_{z}-A_{(F,G)}W-B_{(F,G)}\overline{W}
\label{FGder}%
\end{equation}
if and only if
\begin{equation}
\phi_{\overline{z}}F+\psi_{\overline{z}}G=0. \label{phiFpsiG}%
\end{equation}
This last equation can be rewritten in the following form%
\begin{equation}
W_{\overline{z}}=a_{(F,G)}W+b_{(F,G)}\overline{W} \label{Vekua}%
\end{equation}
which we call the bicomplex Vekua equation. Solutions of this equation are
called $(F,G)$-pseudoanalytic functions.

\begin{remark}
The functions $F$ and $G$ are $(F,G)$-pseudoanalytic, and $\overset{\cdot}%
{F}\equiv\overset{\cdot}{G}\equiv0$.
\end{remark}

\begin{definition}
\label{DefSuccessor}Let $(F,G)$ and $(F_{1},G_{1})$ be two generating pairs in
$\Omega$. $(F_{1},G_{1})$ is called \ successor of $(F,G)$ and $(F,G)$ is
called predecessor of $(F_{1},G_{1})$ if%
\[
a_{(F_{1},G_{1})}=a_{(F,G)}\qquad\text{and}\qquad b_{(F_{1},G_{1})}%
=-B_{(F,G)}\text{.}%
\]

\end{definition}

The importance of this definition becomes obvious from the following statement.

\begin{theorem}
\label{ThBersDer}Let $W$ be an $(F,G)$-pseudoanalytic function and let
$(F_{1},G_{1})$ be a successor of $(F,G)$. Then $\overset{\cdot}{W}$ is an
$(F_{1},G_{1})$-pseudoanalytic function.
\end{theorem}

\begin{definition}
\label{DefAdjoint}Let $(F,G)$ be a generating pair. Its adjoint generating
pair $(F,G)^{\ast}=(F^{\ast},G^{\ast})$ is defined by the formulas%
\[
F^{\ast}=-\frac{2\overline{F}}{F\overline{G}-\overline{F}G},\qquad G^{\ast
}=\frac{2\overline{G}}{F\overline{G}-\overline{F}G}.
\]

\end{definition}

The $(F,G)$-integral is defined as follows
\[
\int_{\Gamma}Wd_{(F,G)}z=F(z_{1})\operatorname{Sc}\int_{\Gamma}G^{\ast
}Wdz+G(z_{1})\operatorname{Sc}\int_{\Gamma}F^{\ast}Wdz
\]
where $\Gamma$ is a rectifiable curve leading from $z_{0}$ to $z_{1}$.

If $W=\phi F+\psi G$ is an $(F,G)$-pseudoanalytic function where $\phi$ and
$\psi$ are scalar functions then
\begin{equation}
\int_{z_{0}}^{z}\overset{\cdot}{W}d_{(F,G)}z=W(z)-\phi(z_{0})F(z)-\psi
(z_{0})G(z), \label{FGAnt}%
\end{equation}
and as $\overset{\cdot}{F}=\overset{}{\overset{\cdot}{G}=}0$, this integral is
path-independent and represents the $(F,G)$-antiderivative of $\overset{\cdot
}{W}$.

\subsection{Generating sequences and Taylor series in formal
powers\label{SubsectGenSeq}}

Following \cite{Berskniga} we introduce the following definitions and results.

\begin{definition}
\label{DefSeq}A sequence of generating pairs $\left\{  (F_{m},G_{m})\right\}
$, $m=0,\pm1,\pm2,\ldots$ is called a generating sequence if $(F_{m+1}%
,G_{m+1})$ is a successor of $(F_{m},G_{m})$. If $(F_{0},G_{0})=(F,G)$, we say
that $(F,G)$ is embedded in $\left\{  (F_{m},G_{m})\right\}  $.
\end{definition}

\begin{theorem}
Let \ $(F,G)$ be a generating pair in $\Omega$. Let $\Omega_{1}$ be a bounded
domain, $\overline{\Omega}_{1}\subset\Omega$. Then $(F,G)$ can be embedded in
a generating sequence in $\Omega_{1}$.
\end{theorem}

\begin{definition}
A generating sequence $\left\{  (F_{m},G_{m})\right\}  $ is said to have
period $\mu>0$ if $(F_{m+\mu},G_{m+\mu})$ is equivalent to $(F_{m},G_{m}),$
that is their characteristic coefficients coincide.
\end{definition}

Let $W$ be an $(F,G)$-pseudoanalytic function. Using a generating sequence in
which $(F,G)$ is embedded we can define the higher derivatives of $W$ by the
recursion formula%
\[
W^{[0]}=W;\qquad W^{[m+1]}=\frac{d_{(F_{m},G_{m})}W^{[m]}}{dz},\quad
m=1,2,\ldots\text{.}%
\]

\begin{definition}
\label{DefFormalPower}The formal power $Z_{m}^{(0)}(a,z_{0};z)$ with center at
$z_{0}\in\Omega$, coefficient $a$ and exponent $0$ is defined as the linear
combination of the generators $F_{m}$, $G_{m}$ with complex constant
coefficients $\lambda$, $\mu$ chosen so that $\lambda F_{m}(z_{0})+\mu
G_{m}(z_{0})=a$. The formal powers with exponents $n=1,2,\ldots$ are defined
by the recursion formula%
\begin{equation}
Z_{m}^{(n+1)}(a,z_{0};z)=(n+1)\int_{z_{0}}^{z}Z_{m+1}^{(n)}(a,z_{0}%
;\zeta)d_{(F_{m},G_{m})}\zeta. \label{recformula}%
\end{equation}

\end{definition}

This definition implies the following properties.

\begin{enumerate}
\item $Z_{m}^{(n)}(a,z_{0};z)$ is an $(F_{m},G_{m})$-pseudoanalytic function
of $z$.

\item If $a^{\prime}$ and $a^{\prime\prime}$ are scalar constants, then
\[
Z_{m}^{(n)}(a^{\prime}+ka^{\prime\prime},z_{0};z)=a^{\prime}Z_{m}%
^{(n)}(1,z_{0};z)+a^{\prime\prime}Z_{m}^{(n)}(k,z_{0};z).
\]

\item The formal powers satisfy the differential relations%
\[
\frac{d_{(F_{m},G_{m})}Z_{m}^{(n)}(a,z_{0};z)}{dz}=nZ_{m+1}^{(n-1)}%
(a,z_{0};z).
\]

\item The asymptotic formulas
\begin{equation}
Z_{m}^{(n)}(a,z_{0};z)\sim a(z-z_{0})^{n},\quad z\rightarrow z_{0}
\label{asymptformulas}%
\end{equation}
hold.
\end{enumerate}

Assume now that
\begin{equation}
W(z)=\sum_{n=0}^{\infty}Z^{(n)}(a,z_{0};z) \label{series}%
\end{equation}
where the absence of the subindex $m$ means that all the formal powers
correspond to the same generating pair $(F,G),$ and the series converges
uniformly in some neighborhood of $z_{0}$. It can be shown that the uniform
limit of pseudoanalytic functions is pseudoanalytic, and that a uniformly
convergent series of $(F,G)$-pseudoanalytic functions can be $(F,G)$%
-differentiated term by term. Hence the function $W$ in (\ref{series}) is
$(F,G)$-pseudoanalytic and its $r$th derivative admits the expansion
\[
W^{[r]}(z)=\sum_{n=r}^{\infty}n(n-1)\cdots(n-r+1)Z_{r}^{(n-r)}(a_{n}%
,z_{0};z).
\]
From this the Taylor formulas for the coefficients are obtained%
\begin{equation}
a_{n}=\frac{W^{[n]}(z_{0})}{n!}. \label{Taylorcoef}%
\end{equation}

\begin{definition}
Let $W(z)$ be a given $(F,G)$-pseudoanalytic function defined for small values
of $\left\vert z-z_{0}\right\vert $. The series%
\begin{equation}
\sum_{n=0}^{\infty}Z^{(n)}(a,z_{0};z) \label{Taylorseries}%
\end{equation}
with the coefficients given by (\ref{Taylorcoef}) is called the Taylor series
of $W$ at $z_{0}$, formed with formal powers.
\end{definition}

The Taylor series always represents the function asymptotically:%
\begin{equation}
W(z)-\sum_{n=0}^{N}Z^{(n)}(a,z_{0};z)=O\left(  \left\vert z-z_{0}\right\vert
^{N+1}\right)  ,\quad z\rightarrow z_{0}, \label{asympt}%
\end{equation}
for all $N$.

If the series (\ref{Taylorseries}) converges uniformly in a neighborhood of
$z_{0}$, it converges to the function $W$.

\subsection{Convergence theorems}

The statements given in this subsection were obtained by L. Bers and S. Agmon
and L. Bers. Their proof in a usual complex case was based on the so called
similarity principle. The similarity principle in general is not valid in a
bicomplex situation. Here we correct the corresponding statement which
unfortunately in \cite{KrAntonio} was formulated with a mistake.

\begin{theorem}
(Similarity principle) Let $w$ be a regular solution of (\ref{Vekua}) in a
domain $\Omega$ such that its values are not zero divisors at any point. Then
the bicomplex function $\Phi=w\cdot e^{h}$, where
\[
h(z)=\frac{1}{\pi}\int_{\Omega}\frac{g(\tau)d\tau}{\tau-z},
\]%
\[
g(z)=\left\{
\begin{array}
[c]{l}%
a_{(F,G)}(z)+b_{(F,G)}(z)\frac{\overline{w}(z)}{w(z)}\text{\quad if }%
w(z)\neq0,\quad z\in\Omega,\\
a_{(F,G)}(z)+b_{(F,G)}(z)\text{\quad if }w(z)=0,\quad z\in\Omega
\end{array}
\right.
\]
is a solution of the equation $\partial_{\overline{z}}\Phi=0$ in $\Omega$.
\end{theorem}

The proof is completely analogous to that for a complex case (see \cite{Vekua}).

Up to now this is an open question how one can guarantee that a solution of
the bicomplex Vekua equation (\ref{Vekua}) be different from a zero divisor at
any point. This is why the proof of the following statements is valid in a
usual complex situation, and in a bicomplex case their validity should be investigated.

In the case when the coefficients in (\ref{Vekua}) are usual complex functions
(with respect to $k$) the following theorems regarding the convergence of
formal Taylor expansions are valid.

\begin{theorem}
\label{ThConvPer}\cite{Berskniga} The formal Taylor expansion
(\ref{Taylorseries}) of a pseudoanalytic function in formal powers defined by
a periodic generating sequence converges in some neighborhood of the center.
\end{theorem}

\begin{definition}
\cite{Berskniga} A generating pair $(F,G)$ is called complete if these
functions are defined and satisfy the H\"{o}lder condition for all finite
values of $z$, the limits $F(\infty)$, $G(\infty)$ exist, $\operatorname{Vec}%
(\overline{F(\infty)}G(\infty))>0$, and the functions $F(1/z)$, $G(1/z)$ also
satisfy the H\"{o}lder condition. \ A complete generating pair is called
normalized if $F(\infty)=1$, $G(\infty)=k$.
\end{definition}

A generating pair equivalent to a complete one is complete, and every complete
generating pair is equivalent to a uniquely determined normalized pair. The
adjoint of a complete (normalized) generating pair is complete (normalized).

From now on we assume that $(F,G)$ is a complete normalized generating pair.
Then much more can be said on the series of corresponding formal powers. We
limit ourselves to the following completeness results (the expansion theorem
and Runge%
\'{}%
s approximation theorem for pseudoanalytic functions).

Following \cite{Berskniga} we shall say that a sequence of functions $W_{n}$
converges normally in a domain $\Omega$ if it converges uniformly on every
bounded closed subdomain of $\Omega$.

\begin{theorem}
\label{ThTaylorRepr}Let $W$ be an $(F,G)$-pseudoanalytic function defined for
$\left\vert z-z_{0}\right\vert <R$. Then it admits a unique expansion of the
form $W(z)=\sum_{n=0}^{\infty}Z^{(n)}(a_{n},z_{0};z)$ which converges normally
for $\left\vert z-z_{0}\right\vert <\theta R$, where $\theta$ is a positive
constant depending on the generating sequence.
\end{theorem}

The first version of this theorem was proved in \cite{AgmonBers}. We follow
here \cite{BersFormalPowers}.

\begin{remark}
\label{RemNecConditions}Necessary and sufficient conditions for the relation
$\theta=1$ are, unfortunately, not known. However, in \cite{BersFormalPowers}
the following sufficient conditions for the case when the generators $(F,G)$
possess partial derivatives are given. One such condition reads:%
\[
\left\vert F_{\overline{z}}(z)\right\vert +\left\vert G_{\overline{z}%
}(z)\right\vert \leq\frac{\operatorname*{Const}}{1+\left\vert z\right\vert
^{1+\varepsilon}}%
\]
for some $\varepsilon>0$. Another condition is
\[
\int\int_{\left\vert z\right\vert <\infty}\left(  \left\vert F_{\overline{z}%
}\right\vert ^{2-\varepsilon}+\left\vert F_{\overline{z}}\right\vert
^{2+\varepsilon}+\left\vert G_{\overline{z}}\right\vert ^{2-\varepsilon
}+\left\vert G_{\overline{z}}\right\vert ^{2+\varepsilon}\right)  dxdy<\infty
\]
for some $0<\varepsilon<1$.
\end{remark}

\begin{theorem}
\cite{BersFormalPowers}\label{ThRunge} A pseudoanalytic function defined in a
simply connected domain can be expanded into a normally convergent series of
formal polynomials (linear combinations of formal powers with positive exponents).
\end{theorem}

\begin{remark}
This theorem admits a direct generalization onto the case of a multiply
connected domain (see \cite{BersFormalPowers}).
\end{remark}

In posterior works \cite{IsmTagieva}, \cite{Menke}, \cite{Fryant} deep results
on interpolation and on the degree of approximation by pseudopolynomials were
obtained. For example,

\begin{theorem}
\cite{Menke}\label{ThMenke} Let $W$ be a pseudoanalytic function in a domain
$\Omega$ bounded by a Jordan curve and satisfy the H\"{o}lder condition on
$\partial\Omega$ with the exponent $\alpha$ ($0<\alpha\leq1$). Then for any
$\varepsilon>0$ and any natural $n$ there exists a pseudopolynomial of order
$n$ satisfying the inequality
\[
\left\vert W(z)-P_{n}(z)\right\vert \leq\frac{\operatorname*{Const}}%
{n^{\alpha-\varepsilon}}\qquad\text{for any }z\in\overline{\Omega}%
\]
where the constant does not depend on $n$, but only on $\varepsilon$.
\end{theorem}

The primary aim of the next section is to show that:

\begin{enumerate}
\item all the mentioned results are of immediate application to the equation
(\ref{maineq}),

\item in many practically important situations the generating sequence and
consequently the formal powers $Z^{(n)}$, $n=0,1,\ldots$ can be constructed explicitly.
\end{enumerate}

\section{Complete systems of solutions for second order equations}

In what follows let us suppose that the scalar function $f$ is defined in a
somewhat bigger domain $\Omega_{\varepsilon}$ with a sufficiently smooth
boundary. Then we change the function $f$ for $z\in\Omega_{\varepsilon}^{{}%
}\backslash\Omega$ and continue it over the whole plane in such a way that
$f\equiv1$ for large $\left\vert z\right\vert $ (see \cite{BersFormalPowers}).
In this way the generating pair $(F,G)=(f,k/f)$ becomes complete and normalized.

Then the following statements are direct corollaries of relations established
in section \ref{SectSchr} between pseudoanalytic functions (solutions of
(\ref{Vekuamain})) and solutions of second order elliptic equations, and
convergence theorems from the previous section.

\begin{definition}
Let $u(z)$ be a given solution of the equation (\ref{maineq}) defined for
small values of $\left\vert z-z_{0}\right\vert $, and let $W(z)$ be a solution
of (\ref{Vekuamain}) constructed according to corollary \ref{CorConjugate}
such that $\operatorname{Sc}W=p^{1/2}u$. The series
\[
p^{-1/2}(z)\sum_{n=0}^{\infty}\operatorname{Sc}Z^{(n)}(a_{n},z_{0};z)
\]
with the coefficients given by (\ref{Taylorcoef}) is called the Taylor series
of $u$ at $z_{0}$, formed with formal powers.
\end{definition}

In the rest of this section we assume that all the coefficients in second
order equations considered in section \ref{SectSchr} are real valued functions
and the particular nonvanishing solution $u_{0}$ of (\ref{maineq}) is real
valued as well.

\begin{theorem}
Let $u(z)$ be a solution of (\ref{maineq}) defined for $\left\vert
z-z_{0}\right\vert <R$. Then it admits a unique expansion of the form
\[
u(z)=p^{-1/2}(z)\sum_{n=0}^{\infty}\operatorname{Sc}Z^{(n)}(a_{n},z_{0};z)
\]
which converges normally for $\left\vert z-z_{0}\right\vert <R$.
\end{theorem}

\begin{proof}
This is a direct consequence of theorem \ref{ThTaylorRepr} and remark
\ref{RemNecConditions}. Both necessary conditions in remark
\ref{RemNecConditions} are fulfilled for the generating pair (\ref{genpair}).
\end{proof}

\begin{theorem}
\label{ThRungeSchr}An arbitrary solution of (\ref{maineq}) defined in a simply
connected domain where there exists a nonvanishing particular solution $u_{0}$
can be expanded into a normally convergent series of formal polynomials
multiplied by $p^{-1/2}$.
\end{theorem}

\begin{proof}
This is a direct corollary of theorem \ref{ThRunge}.
\end{proof}

More precisely the last theorem has the following meaning. Due to Property 2
of formal powers we have that $Z^{(n)}(a,z_{0};z)$ for any Taylor coefficient
$a$ can be easily expressed through $Z^{(n)}(1,z_{0};z)$ and $Z^{(n)}%
(k,z_{0};z)$. Then due to theorem \ref{ThRunge} any solution $W$ of
(\ref{Vekuamain}) can be expanded into a normally convergent series of linear
combinations of $Z^{(n)}(1,z_{0};z)$ and $Z^{(n)}(k,z_{0};z)$. Consequently,
any solution of (\ref{maineq}) can be expanded into a normally convergent
series of linear combinations of scalar parts of $Z^{(n)}(1,z_{0};z)$ and
$Z^{(n)}(k,z_{0};z)$ multiplied by $p^{-1/2}$.

Obviously, for solutions of (\ref{maineq}) the results on the interpolation
and on the degree of approximation like, e.g., theorem \ref{ThMenke} are also valid.

Let us stress that theorem \ref{ThRungeSchr} gives us the following result.
The functions
\begin{equation}
\left\{  p^{-1/2}(z)\operatorname{Sc}Z^{(n)}(1,z_{0};z),\quad p^{-1/2}%
(z)\operatorname{Sc}Z^{(n)}(k,z_{0};z)\right\}  _{n=0}^{\infty}
\label{complsystem}%
\end{equation}
represent a complete system of solutions of (\ref{maineq}) in the sense that
any solution of (\ref{maineq}) can be approximated arbitrarily closely by a
normally convergent series formed by functions (\ref{complsystem}) in any
simply connected domain $\Omega$ where a positive solution of (\ref{maineq})
exists. Moreover, as we show in the next section, in many practically
interesting situations these functions can be constructed explicitly.

\section{Explicit construction of positive formal powers}

The book \cite{Berskniga} (see also \cite[Supplement to Chapter 4]{Courant})
contains explicit formulas for calculation of positive formal powers in the
case when $F$ and $G$ have the form%
\[
F=\left(  \frac{\gamma(x)}{\tau(y)}\right)  ^{1/2}\quad\text{and\quad
}G=k\left(  \frac{\gamma(x)}{\tau(y)}\right)  ^{-1/2}.
\]

In \cite{Krpseudoan} the class of generating pairs for which the generating
sequence and hence the corresponding formal powers can be constructed
explicitly was substantially extended. For the generating pair of the form
(\ref{genpair}) it is possible when $f$ fulfils the following condition.

\begin{condition}
\label{CondRho}(Condition S) \cite{Krpseudoan}\textbf{ }Let $f$ be a scalar
function of some real variable $\rho:$ $f=f(\rho)$ such that the expression
$\frac{\Delta\rho}{\left\vert \nabla\rho\right\vert ^{2}}$ is a function of
$\rho$. We denote it by $s(\rho)=\frac{\Delta\rho}{\left\vert \nabla
\rho\right\vert ^{2}}$.
\end{condition}

Besides the obvious example of any harmonic function $\rho$ and as a
consequence of $\rho$ being a Cartesian variable or $\rho=\arg z=\arctan
(y/x)$, there are many other practically important examples of $\rho$
satisfying Condition S. An important example is $\rho(x,y)=\sqrt{x^{2}+y^{2}}%
$. In this case $s(\rho)=\frac{1}{\rho}$. The parabolic coordinate
$\rho(x,y)=\sqrt{x^{2}+y^{2}}+x$ also fulfills Condition S: $s(\rho)=\frac
{1}{2\rho}$. Elliptic coordinates fulfil Condition S as well (see
\cite{Krpseudoan}).

Denote by $S$ an antiderivative of $s$ with respect to $\rho$.

\begin{theorem}
\cite{Krpseudoan} Let $f$ be a scalar function of a real variable $\rho$
satisfying Condition S and let the function $\varphi=ke^{-S}\rho_{z}$ have no
zeros and be bounded in $\Omega$. Then the generating pair $(F,G)$ with $F=f$
and $G=k/f$ is embedded in the generating sequence $(F_{m},G_{m})$,
$m=0,\pm1,\pm2,\ldots$ with $F=\varphi^{m}F$ and $G=\varphi^{m}G.$
\end{theorem}

This result opens the way for explicit construction of positive formal powers
for the equation (\ref{Vekuamain}) and as a consequence of the complete system
of solutions (\ref{complsystem}) for equation (\ref{maineq}).

Some examples of explicitly constructed formal powers were given in
\cite{Krpseudoan}. Here we show another quite simple but illustrative example.

\begin{example}
Consider the Helmholtz equation
\begin{equation}
(-\Delta+c^{2})u=0 \label{Helmholtz}%
\end{equation}
with $c$ being a real constant. Take the following particular solution of
(\ref{Helmholtz}): $f=e^{cy}$. Let us construct the first few corresponding
formal powers with center at the origin. We have%
\[
Z^{(0)}(1,0;z)=e^{cy},\qquad Z^{(0)}(k,0;z)=ke^{-cy},
\]%
\[
Z^{(1)}(1,0;z)=xe^{cy}+\frac{k\sinh(cy)}{c},\qquad Z^{(1)}(k,0;z)=-\frac
{\sinh(cy)}{c}+kxe^{-cy},
\]

\end{example}%

\[
Z^{(2)}(1,0;z)=\left(  x^{2}-\frac{y}{c}\right)  e^{cy}+\frac{\sinh(cy)}%
{c^{2}}+\frac{2kx\sinh(cy)}{c},
\]%
\[
Z^{(2)}(k,0;z)=-\frac{2x\sinh(cy)}{c}+k\left(  \left(  x^{2}+\frac{y}%
{c}\right)  e^{-cy}-\frac{\sinh(cy)}{c^{2}}\right)  ,\ldots\text{.}%
\]
It is a simple exercise to verify that indeed the asymptotic formulas
(\ref{asymptformulas}) hold. Now taking scalar parts of the formal powers we
obtain a complete system of solutions of the Helmholtz equation:%
\[
u_{1}(x,y)=e^{cy},\qquad u_{2}(x,y)=xe^{cy},\qquad u_{3}(x,y)=-\frac
{\sinh(cy)}{c},
\]%
\[
u_{4}(x,y)=\left(  x^{2}-\frac{y}{c}\right)  e^{cy}+\frac{\sinh(cy)}{c^{2}%
},\qquad u_{5}(x,y)=-\frac{2x\sinh(cy)}{c},\ldots\text{.}%
\]
Formal powers of higher order can be constructed explicitly using a computer
system of symbolic calculation. For this particular example (together with
Maria Rosal\'{\i}a Tenorio) Matlab 6.5 allowed us to obtain analytic
expressions for the formal powers up to the order ten, that gave us the first
twenty one functions $u_{1},\ldots,u_{21}$. We used them for a numerical
solution of the Dirichlet problem for the Helmholtz equation with very
satisfactory results. For example, in the case when $\Omega$ is a unit disk
with centre at the origin, $c=1$ and $u$ on the boundary is equal to $e^{x}$
(this test exact solution gave us the worst precision because of its obvious
\textquotedblleft disparateness\textquotedblright\ from functions $u_{1}%
,u_{2}\ldots$) the maximal error $\max_{z\in\Omega}\left\vert u(z)-\widetilde
{u}(z)\right\vert $ where $u$ is the exact solution and $\widetilde{u}%
=\sum_{n=1}^{21}a_{n}u_{n}$, the real constants $a_{n}$ being found by the
collocation method, was of order 10$^{-7}$. A very fast convergence of the
method was observed.

Although the numerical method based on the usage of explicitly or numerically
constructed pseudoanalytic formal powers still needs a much more detailed
analysis these first results show us that it is quite possible that in due
time and with a further development of symbolic calculation systems it can
rank high among other numerical approaches, especially for solving equations
(\ref{Schrod}) or (\ref{maineq}) with rapidly varying coefficients, when
finite-difference methods fail.

\section{Reduction of the multidimensional second order equation to a first
order equation}

Here we consider the case of dimension $n=3$ and in the final part of this
section we show that a simple generalization gives us the same results in
higher dimensions.

We will consider the algebra $\mathbb{H}(\mathbb{C})$ of complex quaternions
or biquaternions which have the form $Q=Q_{0}+$ $Q_{1}\mathbf{i}%
+Q_{2}\mathbf{j}+Q_{3}\mathbf{k},$ where $\{Q_{k}\}\subset\mathbb{C}$, and
$\mathbf{i}$, $\mathbf{j}$, $\mathbf{k}$ are the quaternionic imaginary units.

The vectorial representation of a complex quaternion will be used. Namely,
each complex quaternion $Q$ is a sum of a scalar $Q_{0}$ and of a vector
\textbf{$Q$}:
\[
Q=\operatorname*{Sc}(Q)+\operatorname*{Vec}(Q)=Q_{0}+\mathbf{Q},
\]
where $\mathbf{Q}=Q_{1}\mathbf{i}+Q_{2}\mathbf{j}+Q_{3}\mathbf{k}$. The
operator of quaternionic conjugation we denote by $C_{H}$: $\overline{Q}%
=C_{H}Q=Q_{0}-\mathbf{Q}$. We conserve the bar for the quaternionic
conjugation which should not provoke any confusion with the same notation for
the conjugation in the first part of the paper because essentially it can be
considered as the same operation if the bicomplex numbers are considered being
embedded in $\mathbb{H}(\mathbb{C})$ in a natural way.

The purely vectorial complex quaternions ($\operatorname*{Sc}(Q)=0$) are
identified with vectors from $\mathbb{C}^{3}$. Note that $\mathbf{Q}%
^{2}=-<\mathbf{Q},\mathbf{Q}>$where $<\mathbf{\cdot},\mathbf{\cdot}>$ denotes
the usual scalar product.

By $M^{P}$ we denote the operator of multiplication by a complex quaternion
$P$ from the right-hand side: $M^{P}Q=Q\cdot P$. More information on the
structure of the algebra of complex quaternions can be found for example in
\cite{AQA} or \cite{KSbook}.

Let $Q$ be a complex quaternion valued differentiable function of
$\mathbf{x}=(x,y,z)$. Denote
\[
DQ=\mathbf{i}\frac{\partial}{\partial x}Q+\mathbf{j}\frac{\partial}{\partial
y}Q+\mathbf{k}\frac{\partial}{\partial z}Q.
\]
This expression can be rewritten in a vector form as follows%

\[
DQ=-\operatorname*{div}\mathbf{Q}+\operatorname*{grad}Q_{0}%
+\operatorname*{rot}\mathbf{Q}.
\]
That is, $\operatorname*{Sc}(DQ)=-\operatorname*{div}\mathbf{Q}$ and
$\operatorname*{Vec}(DQ)=\operatorname*{grad}Q_{0}+\operatorname*{rot}%
\mathbf{Q}$. Let us notice that $D^{2}=-\Delta$. If $Q_{0}$ is a scalar
function then $DQ_{0}$ coincides with $\operatorname*{grad}Q_{0}$.

The following generalization of Leibniz's rule can be proved by a direct
calculation (see \cite[p. 24]{GS1}).

\begin{theorem}
\label{Leib}(Generalized Leibniz rule) Let $\left\{  P,Q\right\}  \subset
C^{1}(G;\mathbb{H}(\mathbb{C}))$, where $G$ is some domain in $\mathbb{R}^{3}%
$. Then
\begin{equation}
D[P\cdot Q]=D[P]\cdot Q+\overline{P}\cdot D[Q]+2(\operatorname{Sc}(PD))[Q],
\label{Leibniz}%
\end{equation}
where
\[
(\operatorname{Sc}(PD))[Q]:=-\sum_{j=1}^{3}P_{j}\partial_{j}Q.
\]

\end{theorem}

We will actively use the following

\begin{remark}
\label{Leibsc}If in Theorem \ref{Leib} $\operatorname{Vec}(P)=0$, that is
$P=P_{0}$, then
\begin{equation}
D[P_{0}\cdot Q]=D[P_{0}]\cdot Q+P_{0}\cdot D[Q]. \label{rem}%
\end{equation}
From this equality we obtain that the operator $D+\frac{\operatorname{grad}%
P_{0}}{P_{0}}$ can be factorized as follows%
\begin{equation}
(D+\frac{\operatorname{grad}P_{0}}{P_{0}})Q=P_{0}^{-1}D(P_{0}\cdot Q).
\label{rem1}%
\end{equation}

\end{remark}

Let $\mathbf{G}$ be a complex valued vector such that $\operatorname{rot}%
\mathbf{G}\equiv0$. Then the complex valued scalar function $\varphi$ is said
to be its antigradient if $\operatorname{grad}\varphi=\mathbf{G}$. We will
write $\varphi=\mathcal{A}[\mathbf{G}]$. The operator $\mathcal{A}$ is a
simple generalization of the usual antiderivative and of the operator
$\overline{A}$ (see section 2), and it defines the function $\varphi$ up to an
arbitrary constant. Its explicit representation is well known and has the form%
\[
\mathcal{A}[\mathbf{G}](x,y,z)=%
{\displaystyle\int\limits_{x_{0}}^{x}}
G_{1}(\xi,y_{0},z_{0})d\xi+%
{\displaystyle\int\limits_{y_{0}}^{y}}
G_{2}(x,\zeta,z_{0})d\zeta+%
{\displaystyle\int\limits_{z_{0}}^{z}}
G_{3}(x,y,\eta)d\eta+C.
\]

Consider the equation%

\begin{equation}
\left(  -\Delta+\nu\right)  g=0\qquad\text{in }G \label{Schr3}%
\end{equation}
where $\Delta=\frac{\partial^{2}}{\partial x^{2}}+\frac{\partial^{2}}{\partial
y^{2}}+\frac{\partial^{2}}{\partial z^{2}}$, $\nu$ and $g$ are complex valued
functions, and $G$ is a domain in $\mathbb{R}^{3}$. We assume that $g$ is
twice continuously differentiable.

\begin{theorem}
\label{Fact3} Let $f$ be a nonvanishing particular solution of (\ref{Schr3}).
Then for any scalar twice continuously differentiable function $g$ the
following equality holds,%
\begin{equation}
(D+M^{\frac{Df}{f}})(D-M^{\frac{Df}{f}})g=\left(  -\Delta+\nu\right)  g.
\label{FactSchr3}%
\end{equation}

\end{theorem}

\begin{proof}
This is a direct calculation based on the Leibniz rule (\ref{rem}).
\end{proof}

\begin{remark}
The factorization (\ref{FactSchr3}) was obtained in \cite{Swansolo},
\cite{Swan} in a form which required a solution of an associated
biquaternionic Riccati equation. In \cite{KKW} it was shown that the solution
has necessarily the form $Df/f$ with $f$ being a solution of (\ref{Schr3}).
\end{remark}

\begin{remark}
Theorem \ref{Fact3} generalizes theorem \ref{Thfact}. In a two-dimensional
situation (\ref{FactSchr3}) reduces to (\ref{fact}).
\end{remark}

\begin{remark}
\label{RemFromSchrToFirst}As $g$ in (\ref{FactSchr3}) is a scalar function,
the factorization of the Schr\"{o}dinger operator can be written in the
following form%
\[
(D+M^{\frac{Df}{f}})fD(f^{-1}g)=\left(  -\Delta+\nu\right)  g,
\]
from which it is obvious that if $g$ is a solution of (\ref{Schr3}) then the
vector $\mathbf{F}=fD(f^{-1}g)$ is a solution of the equation%
\begin{equation}
(D+M^{\frac{Df}{f}})\mathbf{F}=0\qquad\text{in }G. \label{D+MDf}%
\end{equation}
The inverse result we formulate as the following statement.
\end{remark}

\begin{theorem}
\label{ThFromDiracToSchr}Let $\mathbf{F}$ be a solution of (\ref{D+MDf}) in a
simply connected domain $G$. Then $g=f\mathcal{A}[f^{-1}\mathbf{F}]$ is a
solution of (\ref{Schr3}).
\end{theorem}

\begin{proof}
First, in order to apply the operator $\mathcal{A}$ to the vector
$f^{-1}\mathbf{F}$ we should ascertain that indeed,
\begin{equation}
\operatorname*{rot}(f^{-1}\mathbf{F})=0. \label{rot=0}%
\end{equation}
For this, consider the vector part of (\ref{D+MDf}). It has the form
\[
\operatorname*{rot}\mathbf{F}+[\mathbf{F}\times\frac{Df}{f}]=0
\]
which is equivalent to equation (\ref{rot=0}).

Now, applying the Laplacian to $g=f\mathcal{A}[f^{-1}\mathbf{F}]$ and taking
into account that $f$ is a solution of (\ref{Schr3}) and $\mathbf{F}$ is a
solution of (\ref{D+MDf}) we obtain the result:%
\begin{align*}
-\Delta g  &  =D^{2}g=D(Df\cdot\mathcal{A}[f^{-1}\mathbf{F}]+\mathbf{F})\\
&  =f^{-1}\mathbf{F}Df-\mathcal{A}[f^{-1}\mathbf{F}]\Delta f+D\mathbf{F}\\
&  =\mathbf{F}\frac{Df}{f}-\nu f\mathcal{A}[f^{-1}\mathbf{F}]-\mathbf{F}%
\frac{Df}{f}\\
&  =-\nu g.
\end{align*}

\end{proof}

In the same way as in section \ref{SectSchr} we obtain the factorization of
the operator $\operatorname{div}p\operatorname{grad}+q$ where
$\operatorname{div}$ and $\operatorname{grad}$ are already operators with
respect to three independent variables.

\begin{theorem}
Let $u_{0}$ be a nonvanishing particular solution of the equation%
\begin{equation}
(\operatorname{div}p\operatorname{grad}+q)u=0\text{\qquad in }G\subset
\mathbb{R}^{3} \label{maineq3}%
\end{equation}
with $p$, $q$ and $u$ being complex valued functions, $p\in C^{2}(G)$ and
$p\neq0$ in $G$. Then for any scalar function $\varphi\in C^{2}(G)$ the
following equality holds%
\begin{equation}
(\operatorname{div}p\operatorname{grad}+q)\varphi=-p^{1/2}(D+M^{\frac{Df}{f}%
})(D-M^{\frac{Df}{f}})p^{1/2}\varphi\label{mainfact3}%
\end{equation}
where $f=p^{1/2}u_{0}$.
\end{theorem}

\begin{proof}
This is analogous to the proof of theorem \ref{ThFactGenSchr}.
\end{proof}

Thus, if $u$ is a solution of equation (\ref{maineq3}) then
\[
\mathbf{F}=fD(f^{-1}p^{1/2}u)=fD(u_{0}^{-1}u)
\]
is a solution of equation (\ref{D+MDf}) (see remark \ref{RemFromSchrToFirst}).
The inverse result has the following form.

\begin{theorem}
Let $\mathbf{F}$ be a solution of equation (\ref{D+MDf}) in a simply connected
domain $G,$ where $f=p^{1/2}u_{0}$ and $u_{0}$ be a nonvanishing particular
solution of (\ref{maineq3}). Then
\[
u=u_{0}\mathcal{A}[f^{-1}\mathbf{F}]
\]
is a solution of (\ref{maineq3}).
\end{theorem}

\begin{proof}
This is a corollary of theorem \ref{ThFromDiracToSchr} and relation
$(\operatorname{div}p\operatorname{grad}+q)=p^{1/2}(\Delta-\nu)p^{1/2}$ where
$\nu=\Delta f/f$.
\end{proof}

Notice that due to the fact that in (\ref{mainfact3}) $\varphi$ is scalar, we
can rewrite the equality in the form%
\[
(\operatorname{div}p\operatorname{grad}+q)\varphi=-p^{1/2}(D+M^{\frac{Df}{f}%
})(D-\frac{Df}{f}C_{H})p^{1/2}\varphi.
\]
Now, consider the equation
\begin{equation}
(D-\frac{Df}{f}C_{H})W=0, \label{Vekuamain3}%
\end{equation}
where $W$ is an $\mathbb{H}(\mathbb{C})$-valued function. Equation
(\ref{Vekuamain3}) is a direct generalization of the Vekua equation
(\ref{Vekuamain}). Moreover, we show that it preserves some important
properties of (\ref{Vekuamain}).

\begin{theorem}
Let $W=W_{0}+\mathbf{W}$ be a solution of (\ref{Vekuamain3}). Then $W_{0}$ is
a solution of the stationary Schr\"{o}dinger equation
\begin{equation}
-\Delta W_{0}+\nu W_{0}=0, \label{eqSchr3}%
\end{equation}
where $\nu=\Delta f/f$. Moreover, the function $u=f^{-1}W_{0}$ is a solution
of the equation
\begin{equation}
\operatorname{div}(f^{2}\operatorname{grad}u)=0 \label{eqscpart}%
\end{equation}
and the vector function $\mathbf{v}=f\mathbf{W}$ is a solution of the
equation
\begin{equation}
\operatorname{rot}(f^{-2}\operatorname{rot}\mathbf{v})=0. \label{eqvecpart}%
\end{equation}

\end{theorem}

\begin{proof}
Equation (\ref{Vekuamain3}) is equivalent to the system%
\[
\operatorname{div}\mathbf{W}+\left\langle \frac{\nabla f}{f},\mathbf{W}%
\right\rangle =0,
\]%
\[
\operatorname{rot}\mathbf{W}+\left[  \frac{\nabla f}{f}\times\mathbf{W}%
\right]  +\nabla W_{0}-\frac{\nabla f}{f}W_{0}=0
\]
which can be rewritten in the form%
\begin{equation}
\operatorname{div}(f\mathbf{W})=0, \label{syst1}%
\end{equation}%
\begin{equation}
f^{-1}\operatorname{rot}(f\mathbf{W})+f\operatorname{grad}(f^{-1}W_{0})=0.
\label{syst2}%
\end{equation}
From (\ref{syst2}) we obtain (\ref{eqscpart}) and (\ref{eqvecpart}). Equation
(\ref{eqSchr3}) is obtained from (\ref{eqscpart}) and (\ref{factUhlm}).
\end{proof}

\begin{remark}
Observe that the functions
\[
F_{0}=f,\quad F_{1}=\frac{\mathbf{i}}{f},\quad F_{2}=\frac{\mathbf{j}}%
{f},\quad F_{3}=\frac{\mathbf{k}}{f}%
\]
give us a generating quartet for the equation (\ref{Vekuamain3}). They are
solutions of (\ref{Vekuamain3}) and obviously any $\mathbb{H}(\mathbb{C}%
)$-valued function $W$ can be represented in the form%
\[
W=%
{\displaystyle\sum\limits_{j=0}^{3}}
\varphi_{j}F_{j},
\]
where $\varphi_{j}$ are complex valued functions. It is easy to verify that
the function $W$ is a solution of (\ref{Vekuamain3}) iff%
\begin{equation}%
{\displaystyle\sum\limits_{j=0}^{3}}
\left(  D\varphi_{j}\right)  F_{j}=0 \label{Vekuamain3second}%
\end{equation}
in a complete analogy with the two-dimensional case (see remark
\ref{RelAstalaPaivarinta}). Denote
\[
w=\varphi_{0}+\varphi_{1}\mathbf{i}+\varphi_{2}\mathbf{j}+\varphi
_{3}\mathbf{k}.
\]
Then (\ref{Vekuamain3second}) can be written as follows%
\[
D(w+\overline{w})f+D(w-\overline{w})\frac{1}{f}=0
\]
which is equivalent to the equation%
\[
Dw=\frac{1-f^{2}}{1+f^{2}}D\overline{w}.
\]

\end{remark}

\begin{remark}
The results of this section remain valid in the $n$-dimensional situation if
instead of quaternions the Clifford algebra $Cl_{0,n}$ (see, e.g., \cite{BDS},
\cite{GS2}) is considered. The operator $D$ is then introduced as follows
$D=\sum_{j=1}^{n}e_{j}\frac{\partial}{\partial x_{j}}$ where $e_{j}$ are the
basis elements of the Clifford algebra.
\end{remark}

\section{Conclusions}

We showed the possibility of a factorization of the operator
$(\operatorname{div}p\operatorname{grad}+q)$ and investigated only some of its
applications. In a two-dimensional situation under quite general assumptions a
complete system of null solutions of the operator can be constructed
explicitly. It is quite possible that in a multidimensional case using the
results of the preceding section the same can be done. This requires a
multidimensional generalization of L. Bers' theory of formal powers.

Another open question is the proof of expansion and convergence theorems for
the bicomplex Vekua equation of the form (\ref{Vekuamain}).

\end{document}